\documentclass[12pt]{amsart} %这是使用自定义的 aomart 格式模板
\usepackage{geometry}
\usepackage{amsmath}
\usepackage{amssymb}
\usepackage{amsthm}
\usepackage{cite}
\usepackage{graphicx}
\usepackage[footnotesize,centerlast]{subfigure}

\usepackage{enumitem}

\usepackage[all]{xypic}
\usepackage{graphicx}%使用插图的宏包
\usepackage{float} 
\usepackage{color}%使用颜色的宏包
\usepackage{amssymb}%数学符号
\usepackage{amsmath}
\usepackage{yhmath}
\usepackage{mathrsfs}
\usepackage[all]{xy}%交换图
\usepackage{amsthm}
\usepackage{titletoc}
\usepackage{tikz}
\usepackage[colorlinks,linkcolor=blue,anchorcolor=blue,citecolor=blue]{hyperref}% 超链接

\newtheorem{theorem}{Theorem}[section]
\newtheorem{lemma}[theorem]{Lemma}
\newtheorem{proposition}[theorem]{Proposition}
\newtheorem{corollary}[theorem]{Corollary}
\theoremstyle{definition}
\newtheorem{definition}[theorem]{Definition}

\theoremstyle{remark}
\newtheorem{remark}[theorem]{Remark}
\numberwithin{equation}{section}
\newtheorem{clm}{Claim}

\renewcommand\H{{\mathcal H}}

\newcommand{\V}{\operatorname{Vol}}
\newcommand{\M}{\mathcal{M}_g}

\usepackage{microtype}

\begin{document}

\title[Distribution of lengths of closed saddle connections]{Distribution of lengths of closed saddle connections on moduli space of large genus translation surfaces }

\author{Shenxing Zhang}

\address{Shenxing Zhang: School of Mathematics and Statistics, Nanjing University of Science and Technology, Nanjing, 210094, China} \email{52025145@njust.edu.cn}

\date{\today}

\maketitle

%\tableofcontents

\noindent

\begin{abstract}
Let $S_g$ be a closed surface of genus $g$ and $\H_g$ be the moduli space of Abelian differentials on $S_g$.
A stratum of $\H_g$, endowed with the Masur-Veech measure, becomes a probability space. 
Then the number of closed saddle connections with lengths in $[\frac{a}{\sqrt{g}},\frac{b}{\sqrt{g}}]$ on a random translation
surface in the stratum is a random variable.
We prove that when $g\to \infty$, the distribution of the random variable converges to a Poisson distributed random variable. 
This result answers a question of Masur, Rafi and Randecker.

\medskip

\noindent {\bf Keywords:} Random translation surface; closed saddle connection

\medskip

\noindent  {\bf MSC2020:} {30F30, 30F60.}
\end{abstract}

\maketitle

\section{Introduction}
In 2024, Masur, Rafi and Randecker studied the distribution of open saddle connections on a random translation surface in \cite{masur2024lengths}. 
In this paper we study the distribution of closed saddle connections.

A translation surface can be denoted by $(X,\omega)$ where $X$ is a closed Riemann surface of genus $g$ and $\omega$ is an Abelian differential on $X$.
$(X,\omega)$ has a flat structure defined by the conical metric $|\omega|$, with singularities at the zeros of $\omega$.
If an $|\omega|$-geodesic of $(X,\omega)$ connects two zeros (or one zero to itself) of $\omega$ and has no other zeros in its interior,
we call it an open (or closed) saddle connection.

Let $\H_g$ be the moduli space of unit-area translation surfaces of genus $g$.
$\H_g$ has a natural stratification according to the zero order of Abelian differentials.
Particularly, the stratum consisting of all abelian differentials in $\H_g$ which have only simple zeros is called principal stratum, 
denoted by $\H_g(1^{2g-2})$. 

\subsection*{Notation}
For convenience, if a stratum consists of Abelian differentials on which the number of zeros of order $d_i$ is $\kappa_i$
we denote it by 
$$\H_g(d_1^{\kappa_1},\cdots,d_k^{\kappa_k})$$
where $\sum^{k}_{i=1}\kappa_id_i=2g-2$.

Together with a finite natural measure which is introduced by Masur and Veech denoted by $\mu_{MV}$, 
every stratum $\H_g(\kappa)$ of $\H_g$ becomes a probability space
with probability measure $\frac{\mu_{MV}}{\V{\H_g(\kappa)}}$, where $\V{\H_g(\kappa)}$ is the measure of total stratum, 
which we call the Masur-Veech volume of $\H_g(\kappa)$.

Given $(X,\omega)\in \H_g(1^{2g-2})$ and $a,b\in \mathbb{R}^+$,
let $N_{g,[a,b]}(X,\omega)$ be the number of open saddle connections with lengths in the interval $[\frac{a}{g},\frac{b}{g}]$ on $(X,\omega)$.
Then $N_{g,[a,b]}:\H_g(1^{2g-2})\to \mathbb{N}_0$ becomes a random variable.
In \cite{masur2024lengths}, Masur, Rafi and Randecker proved :
\begin{theorem}\label{open}
Let $[a_1,b_1],\cdots,[a_k,b_k]$ be $k$ disjoint intervals. Then as $g\to \infty$, the vector of random variables
$$(N_{g,[a_1,b_1]},\cdots,N_{g,[a_k,b_k]}):\H_g(1^{2g-2})\to \mathbb{N}_0^k$$
converges jointly in distribution to a vector of random variables with Poisson distributions of means $\lambda_{[a_i,b_i]}$, where
$$\lambda_{[a_i,b_i]}=8\pi(b_i^2-a_i^2)$$
for $i=1,\cdots,k$. That is,
$$\mathbb{P}(N_{g,[a_1,b_1]}=n_1,\cdots,N_{g,[a_k,b_k]}=n_k)=\prod_{i=1}^k\frac{\lambda_{[a_i,b_i]}^{n_i}e^{-\lambda_{[a_i,b_i]}}}{n_i!}$$
\end{theorem}

Their method is based on the relationship between Poisson distribution and its factorial moment \cite{bollobas2001random}, 
which was used first in the work of Mirzakhani–Petri \cite{mirzakhani2019lengths} to study the distribution of closed hyperbolic geodesics on random hyperbolic surfaces in Teichm\"uller space with respect to Weil-Petersson measure.

As for translation surface, to use the method, 
\cite{masur2024lengths} introduces a surgery to collapse the open saddle connections to obtain a new general translation surface in a new stratum.
This surgery does not work for closed saddle connections and they proposed the distribution question in the situation of closed saddle connections,
which inspires this work.
This paper gives a new surgery to collapse a closed saddle connection, which ensures the same method can be used for closed saddle connections.

Let $L_{g,[a,b]}(X,\omega)$ be the number of closed saddle connections on $(X,\omega)$ with lengths in $[\frac{a}{\sqrt{g}},\frac{b}{\sqrt{g}}]$.
Denote by $\mathcal{C}_g(\kappa)[\frac{a}{\sqrt{g}},\frac{b}{\sqrt{g}}]$ the subset of $\H_g(\kappa)$ consisting of translation
surfaces on which there exist closed saddle connections
with lengths in $[\frac{a}{\sqrt{g}},\frac{b}{\sqrt{g}}]$.
The main result of this paper is the following theorem 

\begin{theorem}\label{closed}
For the stratum $\H(2^{g-1})$ and disjoint intervals $[a_1,b_1],\cdots,[a_k,b_k]$, when $g\to \infty$, the random variable sequence
$$(L_{g,[a_1,b_1]},\cdots,L_{g,[a_k,b_k]}):\H(2^{g-1}) \to \mathbb{N}^k_0$$
converges jointly in distribution to a vector of random variables with Poisson distributions of means $\lambda_{[a_i,b_i]}$, where
$$\lambda_{[a_i,b_i]}=\frac{3}{2}\pi(b_i^2-a_i^2).$$
\end{theorem}
\begin{remark}
We will explain that the question is trivial for the principal stratum in Section \ref{claim}: 
the probability measure of $\mathcal{C}_g(1^{2g-2})[\frac{a}{\sqrt{g}},\frac{b}{\sqrt{g}}]$ is zero when $g\to \infty$.
For principal stratum, the length interval to be considered should be $[a,b]$.
\end{remark}

Moreover, we consider general stratum $\H(m_1^{O_1(g)},\cdots,m_k^{O_k(g)},1^{O(g)})$, where $\lim\limits_{g\to \infty}\frac{O_i(g)}{g}=c_i$, and we have
\begin{theorem}\label{general m}
Let $\H(m_1^{O_1(g)},\cdots,m_k^{O_k(g)},1^{O(g)})$ be a stratum where $\lim\limits_{g\to \infty}\frac{O_i(g)}{g}=c_i$.
Then for an interval $[a,b]$, $L_{g,[a,b]}(X,\omega)$ converges in distribution to a Poisson distribution of mean $\lambda_{[a,b]}$, where
$$\lambda_{[a,b]}=\sum_{i=1}^{k-1}\frac{c_i}{2}(m_i^2-1)\pi(b^2-a^2).$$
\end{theorem}

\section{Background}
\subsection{Moduli space of Abelian differentials}
Let $S_g$ be a closed surface of genus $g$, $\M$ be the moduli space, i.e, the set of isotopy classes of hyperbolic metric on $S_g$. 

For some $X\in \M$, an Abelian differential $\omega$ on $X$ is a differential $1$-form which under a local coordinate $(U,z)$ of $X$ can be written by
$$\omega=\omega(z)dz,$$
where $\omega(z):U\to \mathbb{C}$ is a holomorphic function. 
The holomorphic differential space on $X$ is denoted by $H^1(X)$, which is a complex linear space of dimension $g$.

The moduli space of Abelian differential $\H_g$ on $S_g$ consists of $(X,\omega)$, 
where $X\in \M$ and $\omega \in H^1(X)$. $\H_g$ is a bundle over $\M$. Since $dimH^1(X)=g, dim\M=3g-3$, we have $dim\H_g=4g-3$.
For $(X,\omega)$, define its area by $\int_X|\omega|^2$.

\subsubsection{Flat metric}

A surface $(X,\omega) \in \H_g$ is endowed with conical metric with singularities on $X$ given by
$$|\omega|=|\omega(z)||dz|.$$
Consider the local coordinate given by $\omega$.
Under the natural coordinate of $\omega$, the transition map is like
$$z\to z+c$$
for a constant $c$. 
We then call $(X,\omega)$ a translation surface and $\H_g$ becomes the moduli space of \emph{translation surfaces}.

Denote the zero set of $(X,\omega)$ by $\Sigma$,
in the natural coordinate, the flat metric becomes $|dz|$ on $X\setminus \Sigma$. 
Around a zero of order $n$, the total angle is $2(n+1)\pi$. Using polar coordinates the metric becomes
$$ds^2=dr^2+(n+1)^2r^2d\theta^2.$$

A \emph{saddle connection} on $(X,\omega)$ is a $|\omega|$- geodesic segment connecting two zeros of $\omega$ and has no other zeros in the interior.
We call a saddle connection \emph{open} if it connects two different zeros and \emph{closed} if it connects a zero to itself.
Generally, a geodesic on $(X,\omega)$ is a concatenation of some saddle connections, with the angle between two adjacent saddle connections larger than $\pi$.
However, there are geodesics that do not pass through any zeros, which we call \emph{regular}.

For a simple closed curve $\gamma$, there are two situations:
\begin{itemize}
  \item[1]  If there are regular geodesics homotopic to $\gamma$, then all such geodesics foliate a flat cylinder, with boundaries not regular.
  \item[2]  There is a unique geodesic homotopic to $\gamma$.
\end{itemize}

For a saddle connection or a closed curve $\gamma$, denote its flat length on $(X,\omega)$ by $\ell_\omega(\gamma)$.
\subsubsection{Stratum and Masur-Veech measure}
For $(X,\omega) \in \H_g$, the sum of its orders of zeros is $2g-2$.
$\H_g$ has a natural stratum structure by the orders of zeros:
let $\kappa=(\kappa_1,\cdots,\kappa_n)$ be a set of $n$ positive integers satisfying $\sum_{i=1}^{n}\kappa_i=2g-2$.
Let $\H_g(\kappa)$ be the subset of $\H_g$ consisting of Abelian differentials with zeros of orders $\kappa_1,\cdots,\kappa_n$.
Then we call $\H_g(\kappa)$ a stratum of $\H_g$ with characterization $\kappa$.
When $\kappa=(1^{2g-2})$, we call $\H_g(\kappa)$ \emph{principal stratum}.
Then the total moduli space $\H_g$ can be decomposed by 
$$\H_g=\bigcup_\kappa \H_g(\kappa),$$
where $\kappa$ takes on all sets of positive integers with sum $2g-2$.

Fix a characterization $\kappa=(\kappa_1,\cdots,\kappa_n)$, the local chart of $\H_g(\kappa)$ can be given by period mapping.
For $(X,\omega)\in \H_g(\kappa)$, denote its zero set by $\Sigma$ and $\{\gamma_1\cdots \gamma_d\}$ be a basis of the relative
homology $H_1(X,\Sigma,\mathbb{Z})$. Then the period mapping can be defined on a neighborhood of $(X,\omega)$ by
$$(X',\omega')\to \left( \int_{\gamma_i}\omega'\right)_{i=1,\cdots,d},$$
where $(X',\omega')$ is in a neighborhood of $(X,\omega)$.
Thus, locally $\H_g(\kappa)$ can be identified with $H_1(X,\Sigma,\mathbb{C})$.

From the dimension of $H_1(X,\Sigma,\mathbb{Z})$, we have
$$dim_{\mathbb{C}}\H_g(\kappa)=dimH_1(X,\Sigma,\mathbb{C})=d=2g-1+n.$$

On each stratum, the period mapping induces a natural Lebesgue measure, normalized so that the lattice $H_1(X,\Sigma,\mathbb{Z}\oplus i\mathbb{Z})$ has co-volume $1$.
This measure, which is introduced by Masur \cite{masur1982interval} and Veech \cite{veech1982gauss}, is called Masur-Veech measure.
For a measured set $K$ of $\H_g(\kappa)$, denote its Masur-Veech measure by $\V(K)$.

Now we consider the hypersurface $\H^1_g(\kappa)$ of $\H_g(\kappa)$ consisting of $(X,\omega)$ with unit area.
The Masur-Veech measure on $\H^1_g(\kappa)$ can be defined as follows:
for an open set $K$ on $\H^1_g(\kappa)$, its Masur-Veech measure is defined by  
$$\V(K)=\operatorname{dim}_{\mathbb{R}}\H_g(\kappa)\cdot \V(Cone(K)),$$
where $Cone(K)$ is the cone of $K$. 
In the following by abuse of notation we still denote $\H^1_g(\kappa)$ by $\H_g(\kappa)$.

\begin{remark}
Note that the stratum $\H_g(\kappa)$ may not be connected, here we consider the connected component. For more details one can refer to \cite{kontsevich2003connected}.
\end{remark}

Aggarwal \cite[Theorem 1.4]{aggarwal2020large} proved that
\begin{theorem}\label{Volume}
For $g>1$ and a stratum $\H_g(\kappa)$ with characterization $\kappa=(\kappa_1,\cdots,\kappa_n)$, 
$$\frac{4}{\prod_{i=1}^n(\kappa_i+1)}\left(1-\frac{2^{2^{200}}}{g} \right)\leq \V(\H_g(\kappa))
\leq \frac{4}{\prod_{i=1}^n(\kappa_i+1)}\left(1+\frac{2^{2^{200}}}{g} \right).$$
\end{theorem}
\subsection{Siegel-Veech constant}
\subsubsection{Configuration}
For a translation surface $(X,\omega)$ and a saddle connection $\gamma$ on it, define the \emph{holonomy vector of $\gamma$} to be $\int_\gamma\omega$ and denote it by $hol(\gamma)$.
Denote the set of all saddle connections of  $(X,\omega)$ by $\tilde{V}_{sc}(X,\omega)$.
Then define
$$V_{sc}(X,\omega)=\{hol(\gamma): \gamma\in\tilde{V}_{sc}(X,\omega)\}.$$
In fact all holonomy vectors of $\gamma \in \tilde{V}_{sc}(X,\omega)$ define a holonomy mapping
$$hol:\tilde{V}_{sc}(X,\omega)\to V_{sc}(X,\omega).$$

By \cite[Proposition 3.1]{vorobets1996planar}, $V_{sc}(X,\omega)$ is a discrete subset of $\mathbb{R}^2$ for any $(X,\omega)$.
Note that different saddle connections may have same holonomy vector.
For example two homologous saddle connections must have same holonomy, 
here \emph{homologous} for two open saddle connections $\beta$ and $\gamma$ means $\beta \circ \gamma^{-1}$ is homologous to zero.
For a vector $\vec{v}\in V_{sc}(X,\omega)$, define the \emph{multiplicity of $\vec{v}$} to be $|hol^{-1}(\vec{v})|$.
It is known that almost any translation surface in any stratum does not have a pair of non-homologous sharing the same holonomy vector 
by \cite[Proposition 3.1]{eskin2003moduli}. So in the following we consider only the surfaces on which if two saddle connections have same holonomy, 
then they are homologous.

To determine how many surfaces have some specific saddle connections, we need to classify the saddle connections by some geometric data,
which we call \emph{configuration}.

For open saddle connections, first fix two zeros $p_1$ and $p_2$ with orders $m_1$ and $m_2$, 
consider the saddle connections connecting $p_1$ and $p_2$ of multiplicity $p$ and denote them by $\gamma_1,\cdots,\gamma_p$,
where the orientation is from $p_1$ to $p_2$ and the cyclic order of $\gamma_i$ is clockwise at $p_1$.
Now we consider the angle between $\gamma_i$ and $\gamma_{i+1}$ at $p_1$ and $p_2$, which are even multiples of $\pi$.
Let the angle between $\gamma_i$ and $\gamma_{i+1}$ at $p_1$ be $2(a'_i+1)\pi$ , 
the angle between $\gamma_i$ and $\gamma_{i+1}$ at $p_2$ be $2(a''_i+1)\pi$.
Then a configuration $\mathcal{C}$ of open saddle connections consists of $(m_1,m_2,a'_i,a''_i)$.

The situation of closed saddle connections is more complicated. Again we consider the closed saddle connections of multiplicity $p$.
Suppose $\gamma_1,\cdots,\gamma_p$ are homologous, and $\gamma_1$ connects $p_1$ to itself.
Note that some $\gamma_i$ may connect some other zero (other than $p_1$) to itself , and some closed saddle connections may bound cylinders.
Denote by $J=\{j_1,\cdots,j_l\}\in \{1,2,\cdots,p\}$ the index set such that $\gamma_{j_1},\cdots, \gamma_{j_l}$ bound $q$ cylinders.
Cut off the $q$ cylinders and $\gamma_1,\cdots,\gamma_p$, there remain $p$ connected components $X_1,\cdots,X_p$. 
$X_1,\cdots,X_p$ can be divided into two types: components of \emph{type one} have two boundary components, components of \emph{type two} have one, 
which is constructed by two $\gamma_i$ and $\gamma_j$ sharing one zero (and don't be cut off).
Consider the angles the boundary curves bound inside the components.
If some $X_k$ is of type one with boundaries $\gamma_k$ and $\gamma_{k+1}$, denote the angles by $(2b'_k+1)\pi$ and $(2b''_k+1)\pi$.
If some $X_i$ is of type two with boundary constructed by $\gamma_i$ and $\gamma_{i+1}$, denote the angles by $2(a'_i+1)\pi$ and $2(a''_i+1)\pi$.
Then a configuration $\mathcal{C}$ of closed saddle connections consists of $(J,b'_k,b''_k,a'_i,a''_i)$.

For a configuration $\mathcal{C}$, denote by $V_{\mathcal{C}}(X,\omega)$ the set of holonomy vectors of saddle connections satisfying $\mathcal{C}$.

\subsubsection{Siegel-Veech formula and Siegel-Veech constant}
Given a stratum $\H_g(\kappa)$ and a configuration $\mathcal{C}$. Let $f:\mathbb{R}^2\to \mathbb{R}$ be an integrable function with compact support, 
the \emph{Siegel-Veech transform} can be defined by
\begin{equation*}
\begin{aligned}
\hat{f}_{\mathcal{C}}&:\H_g(\kappa) \to \mathbb{R} \\
\hat{f}_{\mathcal{C}}(X,\omega)&=\sum_{\vec{v}\in V_{\mathcal{C}}(X,\omega)}f(\vec{v}).
\end{aligned}
\end{equation*}

\begin{theorem}[Siegel-Veech formula \cite{veech1998siegel}]\label{Siegel-Veech}
Let $\H_g(\kappa)$ be a connected stratum and $\mathcal{C}$ a configuration of saddle connections. 
Then there exists a constant $c(\mathcal{C},\H_g(\kappa))$ such that for every integrable 
$f:\mathbb{R}^2\to \mathbb{R}$ with compact support,
$$\mathbb{E}(\hat{f}_{\mathcal{C}})=\frac{1}{\V(\H_g(\kappa))}\int_{\H_g(\kappa)}\hat{f}_{\mathcal{C}}d\mu_{MV}=c(\mathcal{C},\H_g(\kappa))\int_{\mathbb{R}^2}fdxdy.$$

\end{theorem}
The constant $c(\mathcal{C},\H_g(\kappa))$ is called \emph{Siegel-Veech constant} of $\mathcal{C}$ and $\H_g(\kappa)$.

There have been many asymptotics of Siegel–Veech constants for large genus, 
one can see \cite{chen2018quasimodularity,sauvaget2018volumes,aggarwal2019large,aggarwal2020large,chen2020masur} for more results.
Here we only state a result about closed saddle connections, for more details refer to \cite{vallejos2024random}.
\begin{theorem}[\cite{vallejos2024random} Appendix Proposition A.1]\label{Siegel of closed}
Let $\H_g(\kappa)$ be a connected stratum and a fixed zero $\sigma$ of order $m$. 
For any configuration $\mathcal{C}$ of closed saddle connections at $\sigma$ of multiplicity $p$ involving $q$ cylinders and $t \leq p-1$ 
distinct singularities $\sigma_1,\cdots,\sigma_t$ of orders $m_1,\cdots, m_t$ other than $\sigma$, the Siegel-Veech constant associated to $\mathcal{C}$ is 
$$c(\mathcal{C},\H_g(\kappa))=\frac{(m+1)\prod_{i=1}^t(m_i+1)}{2^{p-1}|\Gamma_{-}||\Gamma|(p+q-1)!}(1+o(1))^pO(\frac{1}{g^{2p+q-2}}),$$
where $|\Gamma_{-}|$ and $|\Gamma|$ are combinatorial constants depending on the combinatorial ways of these cylinders.
\end{theorem}

\subsection{Factorial moment}
For a probability space $(\Omega,\mathbb{P})$ and a random variable $X$ on it, its factorial moment generating function (FMGF) is defined as:
$$M_X(t)=\mathbb{E}[(1+t)^X],$$
and its k-th factorial moment $\mathbb{E}[(X)_k]$ is defined by
$$\mathbb{E}[(X)_k]=\mathbb{E}[X(X-1)\cdots(X-k+1)].$$
In fact it is the k-th derivative of $M_X(t)$ with respect to $t$ at $t=0$.

Recall that a random variable $X:\Omega \to \mathbb{N}_0$ is Poisson distributed with mean $\lambda$ if
$$\mathbb{P}(X=k)=\frac{\lambda^k e^{-\lambda}}{k!}.$$
The following theorem from \cite{bollobas2001random} gives the relationship between factorial moment and distribution for Poisson distribution:
\begin{theorem}[The method of moment from \cite{bollobas2001random}, Theorem 1.23]\label{moment}
Let $\left\{(\Omega_i,\mathbb{P}_i)\right\}_{i\in \mathbb{N}}$ be a sequence of probability spaces. 
For $m\in \mathbb{N}$, let $N_{1,i},\cdots,N_{m,i}:\Omega_i\to \mathbb{N}_0$ be random variables for all $i\in \mathbb{N}$
and suppose there exist $\lambda_1,\cdots,\lambda_m\in(0,\infty)$ such that:
$$\lim_{i\to \infty}\mathbb{E}[(N_{1,i})_{k_1}\cdots(N_{m,i})_{k_m}]=\lambda_1^{k_1}\cdots\lambda_m^{k_m}$$
for all $k_1,\cdots,k_m\in \mathbb{N}$. Then
$$\lim_{i\to \infty}\mathbb{P}[(N_{1,i})=n_1,\cdots,(N_{m,i})=n_m]=\prod^k_{j=1}\frac{\lambda_j^{n_j} e^{-\lambda_j}}{n_j!}.$$
That is, $(N_{1,i},\cdots,N_{m,i}):\Omega_i \to \mathbb{N}^m_0$ converges jointly in distribution to a vector which is independently Poisson distributed with means 
$\lambda_1,\cdots,\lambda_m$.

\end{theorem}

\section{Some claims}\label{claim}
Before the proof of the main theorem, we have to make some claims which will simplify the situation.
We consider the stratum whose zeros have order of $O(1)$, and it implies $|\kappa|=O(g)$.

\begin{clm} \label{clm1}
A closed saddle connection is non-separating.
\end{clm}
\begin{proof}
For a closed saddle connection $\gamma$ of $(X,\omega)$, if it's separable, then it cuts $X$ into $X_1$ and $X_2$, and $\partial X_1=\partial X_2=\gamma$.
By Stokes' theorem, $\int_\gamma \omega=0$, but since $\gamma$ is a saddle connection its length is $|\int_\gamma \omega|$, which is a contradiction.
\end{proof}

\begin{clm} \label{clm2}
The probability measure of the subset of $\H_g(\kappa)$ which has a closed saddle connection involving $q$ cylinders of multiplicity $p$ 
goes to $0$ as $g\to \infty$, where $p>1$ or $q>0$.
\end{clm}

\begin{proof}
From Theorem \ref{Siegel of closed}, 
the Siegel-Veech constant of configuration with closed saddle connection involving $q$ cylinders of multiplicity $p$ is $O(\frac{1}{g^{2p+q-2}})$ . 
Denote by $\mathcal{C}^{p,q}_g(\kappa)[\frac{a}{\sqrt{g}},\frac{b}{\sqrt{g}}]$ the subset of $\H_g(\kappa)$ consisting of surfaces which have a closed saddle connection 
of length in the interval $[\frac{a}{\sqrt{g}},\frac{b}{\sqrt{g}}]$
involving $q$ cylinders of multiplicity $p$.
Then for $p\geq 2$ or $q\geq 1$, we have
 $$\V(\mathcal{C}^{p,q}_g(\kappa)[\frac{a}{\sqrt{g}},\frac{b}{\sqrt{g}}])\leq |\kappa|O(\frac{1}{g^{q}})\frac{b^2-a^2}{g}\V(\H_g(\kappa))=O(\frac{1}{g})\V(\H_g(\kappa)$$
where the coefficient $|\kappa|$ is the selection of the zero that the closed saddle connection goes through. 
So when $g\to \infty$ its probability measure goes to $0$ 
and it suffices to consider the closed saddle connection with multiplicity $1$ and no cylinders around it.
\end{proof}

\begin{clm} \label{clm3}
The angle of a closed saddle connection at the zero which it connects is an odd multiple of $\pi$.
\end{clm}

\begin{proof}
For an Abelian differential $(X,\omega)$, if it has a closed saddle connection $\gamma$ connecting a zero $p$, if its angle at $p$ is $2k\pi$,
then the holonomy of $\gamma$ and $-\gamma$ have same direction. But
$$\int_\gamma \omega=-\int_{-\gamma} \omega.$$
It implies $\int_\gamma \omega=0$, which is a contradiction.
\end{proof}

\begin{clm} \label{clm4}
If a closed saddle connection has angle $\pi$ at one side, it must have a cylinder at the side.
\end{clm}

\begin{proof}
Let $\gamma$ be a colsed saddle connection on $(X,\omega)$ which has angle $\pi$ at one side. Since $(X,\omega)$ is oriented, its normal bundle of $\gamma$ at the side is oriented.
Consider the exponential map from the bundle to $(X,\omega)$, which is well-defined on $\gamma \times [0,s]$ for some $s$ sufficiently small and the image is a cylinder.

\end{proof}

Combining claim $2,3$ and $4$, we have that for a closed saddle connection $\gamma$ on $(X,\omega)$ in principal stratum, 
since its total angle around the zero is $4\pi$ and has to be divided into odd multiples of $\pi$, 
it must have angle $\pi$ at one side of $\gamma$ and $3\pi$ at another.
So there is a cyliner with $\gamma$ as a boundary. As for the other boundary, if it consists of some open saddle connections, 
then these saddle connections have angles $\pi$. That means there are at least two non-homologous saddle connections have angle $\pi$, locally from the period coordinate, the surfaces in this situation consist of a hyperplane of the stratum
so is a set of measure zero. 
So we only need to consider the cylinder bounded by curves homologous to $\gamma$.
Then from Claim \ref{clm2} the probability measure of $\mathcal{C}_g(1^{2g-2})[\frac{a}{\sqrt{g}},\frac{b}{\sqrt{g}}]$ has limit
 $$\lim_{g\to \infty}\frac{\V(\mathcal{C}_g(1^{2g-2})[\frac{a}{\sqrt{g}},\frac{b}{\sqrt{g}}])}{\V(\H_g(1^{2g-2})}=O(\frac{1}{g})\to 0.$$

So we need to consider the question on stratum with higher order zeros. Next we will consider $\H_g(2^{g-1})$ first.

\section{Surgery}\label{surgery}
In this section we will introduce some surgeries that, when combined, can collapse a closed saddle connection.
Then we will use a series of such surgeries to construct a mapping to collapse a closed saddle connection.
\subsection{Open up a higher-order zero}
Masur, Rafi and Randecker introduce a surgery in \cite{masur2024lengths} which is a variation of a surgery from \cite{eskin2003moduli}. 
This surgery can collapse a saddle connection which is not closed and has multiplicity $1$.

Let $(X,\omega)$ be an Abelian differential, $\sigma$ be an open saddle connection with endpoints $v_1$ and $v_2$, 
whose orders are $n_1$ and $n_2$ respectively.
Since the total angle at $v_1$ is $2(n_1+1)\pi$, 
one can extend $\sigma$ from $v_1$ along its direction to $v_3$ and denote by $\sigma'$ the geodesic segment $v_1v_3$ 
such that $\ell_\omega(\sigma)=\ell_\omega(\sigma')$
and the angle between $\sigma$ and $\sigma'$ is
$2k_1\pi$ and $2k_2\pi$, where $k_1+k_2=n_1+1$. 
If $\sigma'$ does not go through other zeros of $(X,\omega)$, 
which means $\sigma'$ is a ray from $v_1$, the following surgery can be carried out.

Cut along $\sigma+\sigma'$ and denote the two copies of $\sigma+\sigma'$ and $v_1$ by $\sigma^{\pm}+(\sigma')^{\pm}$ and $v_1^{\pm}$. 
Then glue $\sigma^+$ and $(\sigma')^+$, $\sigma^-$ and $(\sigma')^-$. 
This surgery reduces the order of $v_1$ and constructs a new Abelian differential $(X',\omega')$. 
On $(X',\omega')$, the total angles at $v_1^+$, $v_1^-$ and $v_2$ are $2k_1\pi$, $2k_2\pi$, and $2(n_2+2)\pi$ respectively,
that is, the orders of $v_1^+$, $v_1^-$ and $v_2$ are $k_1-1$, $k_2-1$, and $n_2+1$.
Particularly, if $v_1$ is a simple zero, the surgery collapses the saddle connection. We call this surgery \emph{collapsing surgery}.
\begin{figure}[h]
    \centering
    \includegraphics[width=0.6\linewidth]{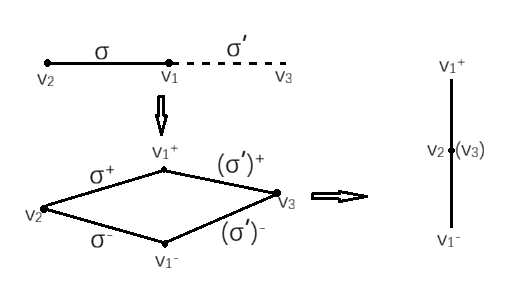}
    \caption{The collapsing surgery}
    \label{collapsing}
\end{figure}

In our situation, we need to use the inverse of collapsing surgery which we call \emph{opening surgery}.
From Claim \ref{clm2}, it suffices to consider the subset $\mathcal{C}^{1,0}_g(2^{g-1})[\frac{a}{\sqrt{g}},\frac{b}{\sqrt{g}}]$.
For an Abelian differential $(X,\omega)$ in $\mathcal{C}^{1,0}_g(2^{g-1})[\frac{a}{\sqrt{g}},\frac{b}{\sqrt{g}}]$, 
choose a closed saddle connection $\gamma$ with length in $[\frac{a}{\sqrt{g}},\frac{b}{\sqrt{g}}]$ at a double zero $p$.
We can reverse the surgery introduced above as follows.

Since the total angle at $p$ is $6\pi$, from Claim \ref{clm3} and Claim \ref{clm4}, the angles on both sides of $\gamma$ at $p$ must be $3\pi$.
Fix the orientation of $\gamma$, choose two rays $\sigma^+$ and $\sigma^-$ from $p$ such that 
$\ell_\omega(\sigma^+)=\ell_\omega(\sigma^-)=\ell_\omega(\gamma)$,
and the angle between $\gamma$ and $\sigma^+$
is $\pi$, the angle between $\gamma$ and $\sigma^-$ is $-\pi$, where the sign is consistent with the orientation of the surface.
Denote by $q^+$ and $q^-$ the other endpoints of $\sigma^+$ and $\sigma^-$.
Then cut along $\sigma^+ + \sigma^-$ and glue the two copies of $\sigma^+$ and $\sigma^-$.
After the surgery, the two copies of $p$ become a regular and a simple zero, $q^+$ and $q^-$ are glued to become a new simple zero $q$.

\begin{figure}[h]
   \centering
   \includegraphics[width=0.6\linewidth]{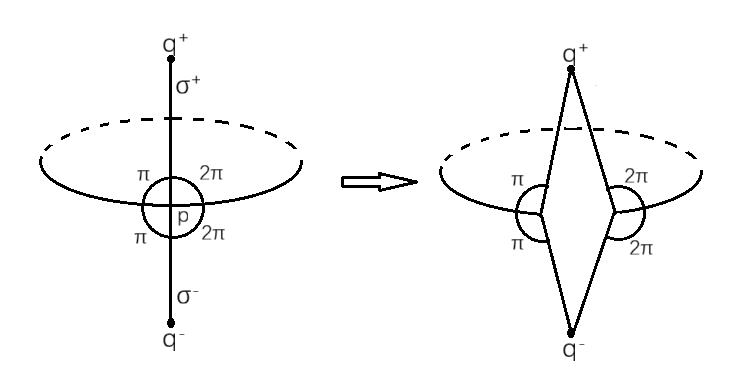}
  \caption{The opening surgery}
  \label{opening}
\end{figure}

This surgery is the inverse of the collapsing surgery above: the double zero is replaced by two simple zeros,
and it can be carried out except for the surgery locus going through some zero, 
that is, it has two non-homologous saddle connections with angle $2\pi$.
From period mapping, such a subset has measure zero. We call the subset on which the opening surgery can be carried out \emph{permissible set}, 
and denote it by $\mathcal{C}^{1,0}_{g,OP}(2^{g-1})[\frac{a}{\sqrt{g}},\frac{b}{\sqrt{g}}]$, from the discussion above we have 
\begin{proposition}\label{permission}
$\mathcal{C}^{1,0}_{g,OP}(2^{g-1})[\frac{a}{\sqrt{g}},\frac{b}{\sqrt{g}}]$ is a full measure subset of $\mathcal{C}^{1,0}_g(2^{g-1})[\frac{a}{\sqrt{g}},\frac{b}{\sqrt{g}}]$.
\end{proposition}

%By the surgery of collapsing this map is one-to-one and denote the opening surgery by $F_1$.
\begin{figure}[h]
    \centering
    \includegraphics[width=0.6\linewidth]{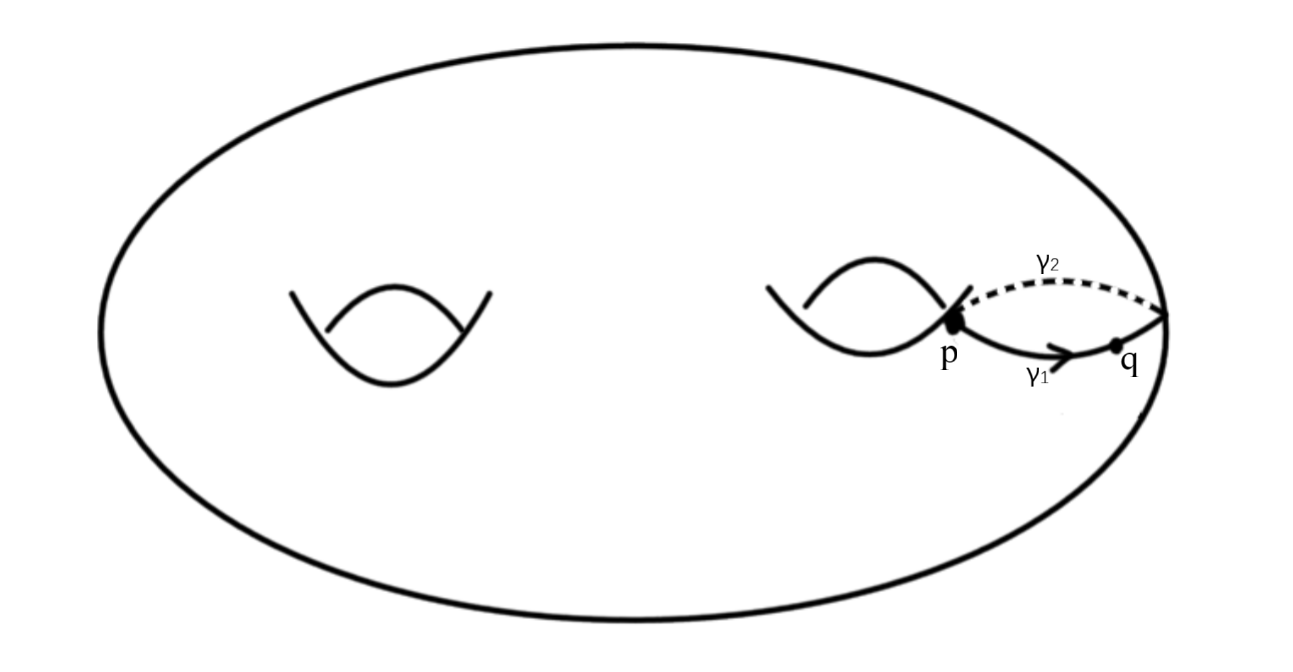}
    \caption{The resulting surface}
    \label{opening surface}
\end{figure}

\begin{remark}
We can also choose $\sigma^{\pm}$ to make the angle between $\gamma$ and them are $2\pi$, 
it suffices to make sure the opening operation can obtain an abelian differential and a smooth loop homotopic to $\gamma$.
\end{remark}
\subsection{Move zero along closed curve and pinch}
%For some $(X,\omega)\in \mathcal{C}^{1,0}_{g,Per}(2^{g-1})$, let $(X',\omega')=F_1(X,\omega)$.
We have constructed the opening surgery, 
which takes a closed saddle connection $\gamma$ at a double zero to two saddle connections $\gamma_1$ and $\gamma_2$
sharing the same endpoints $p$ and $q$ which are both simple zeros.
Moreover, the angles between the two saddle connections at the two simple zeros (on both sides) are all $2\pi$, 
and if we denote by $\ell_\omega(\gamma)=L$, we have $\ell_{\omega'}(\gamma_1)=L$, $\ell_{\omega'}(\gamma_2)=2L$, where $(X',\omega')$ is the resulting surface of $(X,\omega)$ after opening surgery along $\gamma$.
If we want to collapse $\gamma_1$ and $\gamma_2$ simultaneously, we have to make a surgery to move $q$ along $\gamma_2$ to adjust the length of $\gamma_2$.

\subsubsection{Move the zero locally by period mapping.}
Let $(X',\omega')$ be the resulting surface of $(X,\omega)$ after opening surgery for some closed saddle connection $\gamma$.
First we need to choose a special basis of $H_1(X',\Sigma;\mathbb{C})$, where $\Sigma=\{p,q,p_2,\cdots,p_{g-2}\}$
is the zero set of $(X',\omega')$.
From Claim \ref{clm1}, the homology class $[\gamma]=[\gamma_1+\gamma_2]$ is non-separating.
Let $\tilde{\alpha}_1=[\gamma]$.

We will choose the basis for the marked Abelian differentials.
Lift $\H_g$ to a fundamental domain $\tilde{\H_g}$ on its universal covering, in which an
element is like $(X', f: S\to X', \omega')$, where $(X', f: S\to X')$ is a marked Riemann surface and $\omega'$ is an abelian differential on $X'$, 
$S$ is a closed topological surface of genus $g$. Choose a canonical basis $a_1,b_1,\cdots,a_{g},b_{g}$ on $S$, where $a_i$ is dual to $b_i$, and the intersection number becomes a symplectic form $w(\cdot,\cdot)$ on $H_1(S;\mathbb{Z})$.

We need some algorithms to find our basis, for details on can refer to \cite{Cohen1993}.
Suppose $\alpha_1=f^{-1}(\tilde{\alpha}_1)$ has coordinate $(x_1,y_1,\cdots,x_{g},y_g)$ under the canonical basis $\{a_1,b_1,\cdots,a_{g},b_{g}\}$ in $H_1(S;\mathbb{Z})$. Consider the following equation in terms of 
$\beta_1=(x'_1,y'_1,\cdots,x'_{g},y'_g)$:
$$w(\alpha_1,\beta_1)=1.$$
Use the Euclid's Extended Algorithm in \cite[$\S$1.3.2]{Cohen1993}, we can find a primitive solution of
$\beta_1 \in H_1(S;\mathbb{Z})$, which represents a curve dual to $\alpha_1$. Let $V$ be the subspace generated by $\alpha_1$ and $\beta_1$.
Consider the symplectic complement $V'$ of $V$ under $w$ in $H_1(S;\mathbb{Z})$, this is
an integer lattice. Choose a basis of $V'$ to construct a matrix $A$, and use the Hermite
Normal Form Algorithm \cite[$\S$2.4.2]{Cohen1993} to find the unique Hermite normal basis
$(\hat{\alpha_2},\hat{\beta_2},\cdots,\hat{\alpha_g},\hat{\beta_g})$ of $V'$, such that
the first basis $\hat{\alpha_2}$ is primitive.

Note that $(\hat{\alpha_2},\hat{\beta_2},\cdots,\hat{\alpha_g},\hat{\beta_g})$ may not become a symplectic basis of $V'$. Let $\alpha_2=\hat{\alpha_2}$, and
we use Euclid's Extended Algorithm to find a new vector $\beta_2$ in $V'$
satisfying $w(\alpha_2,\beta_2)=1$.
Next consider the subspace  generated by $\alpha_2$ and $\beta_2$ and its
symplectic complement $V''$,
repeat this process we extend a symplectic basis of $H_1(S;\mathbb{Z})$
$$(\alpha_1,\beta_1,\cdots,\alpha_g,\beta_g)$$
such that $f(\alpha_1)=\tilde{\alpha_1}$, and $f$ maps the basis to a canonical basis
of $H_1(X';\mathbb{Z})$, which we still denote by
$$(\alpha_1,\beta_1,\cdots,\alpha_g,\beta_g).$$

\begin{remark}
In fact by this way we can find a canonical basis for each surface in the interior of $\tilde{\H}_g$.
\end{remark}

In $\tilde{\H_g}$ choose the fundamental domain $\tilde{\H_g}(2^{g-1})$ of $\H_g(2^{g-1})$, and 
lift the set $\mathcal{C}^{1,0}_{g,OP}(2^{g-1})[\frac{a}{\sqrt{g}},\frac{b}{\sqrt{g}}]$ to $\tilde{\mathcal{C}}^{1,0}_{g,OP}(2^{g-1})[\frac{a}{\sqrt{g}},\frac{b}{\sqrt{g}}]$
on $\tilde{\H_g}(\kappa)$.
The basis we have chosen is well-defined on the interior of $\tilde{\mathcal{C}}^{1,0}_{g,Per}(2^{g-1})[\frac{a}{\sqrt{g}},\frac{b}{\sqrt{g}}]$,
which projects to a full-measure subset of $\mathcal{C}^{1,0}_{g,Per}(2^{g-1})[\frac{a}{\sqrt{g}},\frac{b}{\sqrt{g}}]$.
We still denote the subset by $\mathcal{C}^{1,0}_{g,Per}(2^{g-1})[\frac{a}{\sqrt{g}},\frac{b}{\sqrt{g}}]$.
Finally we get the new permissible set by removing the projection of $\partial{\tilde{\H_g}(\kappa)}$ from the old permissible set, 
so for every surface in the permissible set we can choose such basis .

Next we choose the relative homology class
$$(pq,pp_2\cdots,pp_{g-2}),$$
where $[\gamma_1]=pq$, and $pp_i$ is  homotopic to $0$ as free homotopy class in $H_1(X',\mathbb{C})$. 
Together they compose a basis of $H_1(X',\Sigma;\mathbb{C})$, 
and their holonomy gives a period mapping around $(X',\omega')$ locally.

Let $\omega'(\alpha_1)=(x_1,y_1)$, $\omega'(\beta_1)=(x_2,y_2)$, $\omega'(\alpha_2)=(x_3,y_3)$,
$\omega'(\beta_2)=(x_4,y_4)$ and suppose $x_1y_2-x_2y_1>0, y_4>0$.
Now choose a relative cohomology class $\upsilon \in H^1(X,\Sigma;\mathbb{C})$ ,
which can be considered as a tangent vector in $T_{\tilde{\H_g}(1,1,\cdots,2)}(X',f, \omega')$ such that
\begin{equation}\label{moving}
\upsilon(\alpha_1)=(-x_1,-y_1), \upsilon(\alpha_2)=(\frac{y_2x_1}{y_4}-\frac{x_2y_1}{y_4},0)
\end{equation}
and on the other basis we assign zero to $\upsilon$. 
On the interior of $\tilde{\H_g}$, locally from period mapping $\upsilon$ is a well-defined vector field,
so it can be lifted to the permissible set to become a well-defined vector field.

\begin{remark}\label{explain}
Here $\upsilon(\alpha_1)$ is chosen to guarantee that along the curve $(X',\omega')+t\upsilon$ in moduli space the zero moves along $\gamma$, 
and $\upsilon(\alpha_2)$ is chosen to guarantee the resulting surface has area $1$, since the area of $(X,\omega)$ can be written by
$$\int_X|\omega|^2=\frac{i}{2}\sum_i(\int_{\alpha_i}\omega \int_{\beta_i}\overline{\omega}-\int_{\beta_i}\omega \int_{\alpha_i}\overline{\omega}).$$
\end{remark}
%which gives 

Consider the curve in $\H_g(1,1,2,\cdots,2)$: $(X', \omega')+t\upsilon, t\in [0,1]$. 
If $\gamma_1$ and $\gamma_2$ don't degenerate along the curve, then we obtain a new surface $(X',\omega')+\upsilon$.
\iffalse
Let $\alpha_1=[\gamma_1+\gamma_2]$ and  $\beta_2$ be the curve whose flat geodesic has minimal vertical projection during 
$(M,f)$.
Since the vertical projection of all saddle connections is not changed along $(X',\omega')+t\upsilon$, 
the holonomy class of $\beta_2$ is unchanged.

So when a closed saddle connection $\gamma$ is chosen, at any time $t_0$ we can always choose $\beta_2$ and $\upsilon$, which guarantees\fi
Since the homology class of the basis is not changed along $(X',\omega')+\upsilon$,
we can choose the same basis at any time $t$ on $(X',\omega')+\upsilon$, which means
$\upsilon$ is indepedent to the locally period mapping becomes a vector field globally.
%$$F_2:\mathcal{C}^{1,0}_{g,Per}(2^{g-1})[\frac{a}{\sqrt{g}},\frac{b}{\sqrt{g}}]  \to \H_g(1^2,2^{g-1})$$
%such that $F_2[(X,\omega)]=F_1[(X,\omega)]+\upsilon$. 
From the construction above the resulting surface has two saddle connections 
in the relative homology class $\gamma_1$ and $\gamma_2$ with equal length.

\subsubsection{Moving surgery is well-defined almost everywhere}\label{F_2}
Let $(X,\omega)$ be a surface in $ \mathcal{C}^{1,0}_{g,Per}(2^{g-1})[\frac{a}{\sqrt{g}},\frac{b}{\sqrt{g}}]$.
After opening a closed saddle connection $\gamma$, $(X,\omega)$ becomes $(X',\omega')$.
We have defined the moving surgery if the saddle connections $\gamma_1$ and $\gamma_2$ are preserved along the curve $(X',\omega')+t\upsilon, t\in [0,1]$.
Next we will see this surgery can be defined for almost every translation surface in the set of resulting surfaces under opening surgery on $\mathcal{C}^{1,0}_{g,Per}(2^{g-1})$.

From above we know the moving surgery can be realized except when $\gamma_1$ or $\gamma_2$ degenerates along $(X',\omega')+t\upsilon, t\in [0,1]$.
Suppose $T$ is the first time when the geodesic in $[\gamma_1]$ can be represented by $pp_i-p_ip_j-\cdots-p_kq$,
since the saddle connection is smooth on $(X',\omega')+t\upsilon,t\in [0,T]$, 
we must have 
$$\vec{pp_i}+\vec{p_ip_j}+\cdots+\vec{p_kq}=\omega'+T\upsilon(\gamma_1),$$ 
$$|pp_i|+\cdots+|p_kq|=|\gamma_1|.$$
So the corner of $\vec{pp_i}-\vec{p_ip_j}-\cdots-\vec{p_kq}$ must be $\pi$. 
But the holonomy of $pp_i$ is not changed on $(X',\omega')+t\upsilon,t\in [0,T]$
This implies there exists some $pp_i$ such that its holonomy on $(X',\omega')$ has the same direction with $p_0p_1$. 
And by the relation of holonomy under opening surgery, this implies $p_0p_i$ has the same direction with the closed curve on $(X,\omega)$.
Under period mapping, on each local chart this is a measure-zero set.
So all such $(X,\omega)$ are a subset of measure zero in $ \H_g(2^{g-1})$, also in $\mathcal{C}^{1,0}_{g,Per}(2^{g-1})[\frac{a}{\sqrt{g}},\frac{b}{\sqrt{g}}]$.

Moreover, if some $(Y,\theta)$ is a resulting surface of moving surgery, note that there exists at least one simple closed curve consisting of a pair of non-homologous saddle
connections sharing same endpoints, which we denote by $\alpha_1$.
Choose the preimage $(Y,f,\theta)$ of $(Y,\theta)$ in $\tilde{\H_g}$ under the covering.
By the construction of moving surgery, we can also find the canonical basis on $Y$ construct the tangent vector $\upsilon$, where $\alpha_1$ is a homology class consists of two non-homologous saddle connections which has same endpoints and same lengths, 
since $(Y,\theta)$ is in  the image, there exists at least one pair of such saddle connections.
And from the construction if $(Y,\theta)$ has two preimages it must have at least two such pairs of saddle connections.

\subsubsection{Pinching}
On the subset of $\mathcal{C}^{1,0}_{g,Per}(2^{g-1})[\frac{a}{\sqrt{g}},\frac{b}{\sqrt{g}}]$ where the moving surgery is well-defined,
we can cut along $\gamma_1+\gamma_2$ and denote the two copies of $p$ and $q$ by $q^+$, $q^-$, $p^+$, $p^-$, 
the two copies of $\gamma_1$ and $\gamma_2$ by $\gamma^+_1$, $\gamma^-_1$, $\gamma^+_2$, $\gamma^-_2$.
Notice that the total angles at $q^+$, $q^-$, $p^+$, $p^-$ are all $2\pi$, 
which means $\gamma^\pm_1$ have same direction with $\gamma^\pm_2$,
then glue $\gamma^+_1$ with $\gamma^+_2$, $\gamma^-_1$ with $\gamma^-_2$, which make $q^+$, $q^-$, $p^+$, $p^-$ become regular points.
Since $\gamma_1+\gamma_2$ is non-separating, this surgery reduces the genus by one, and obtains a new abelian differential 
on which we mark two points $p^+$ and $p^-$ so as to find its inverse. 
%Considering the choice of which zero is selected after applying the inverse map, we finally obtain:
%$$F_3: \mathcal{C}^{1,0}_{g,Per}(2^{g-1})[\frac{a}{\sqrt{g}},\frac{b}{\sqrt{g}}]\to M\times \H_{g-1}(2^{g-2},0,0)\times A_{[\frac{a}{\sqrt{g}},\frac{b}{\sqrt{g}}]},$$
%where $M$ is the combinatorial data of choosing one from $g-1$ zeros.

\begin{figure}[h]
    \centering
    \includegraphics[width=0.6\linewidth]{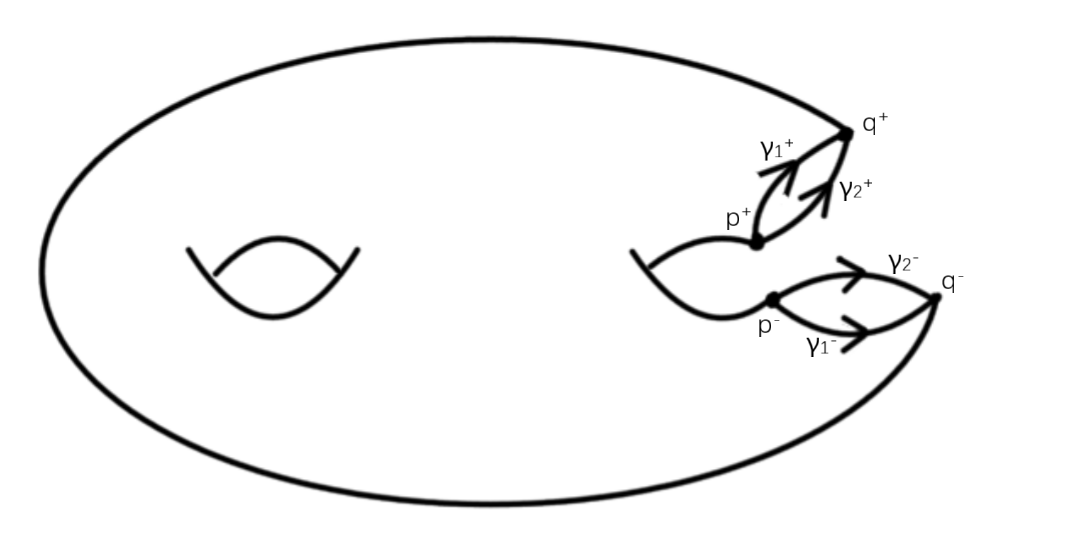}
    \caption{Moving and cut}
    \label{moving and cut}
\end{figure}

\begin{figure}[h]
    \centering
    \includegraphics[width=0.6\linewidth]{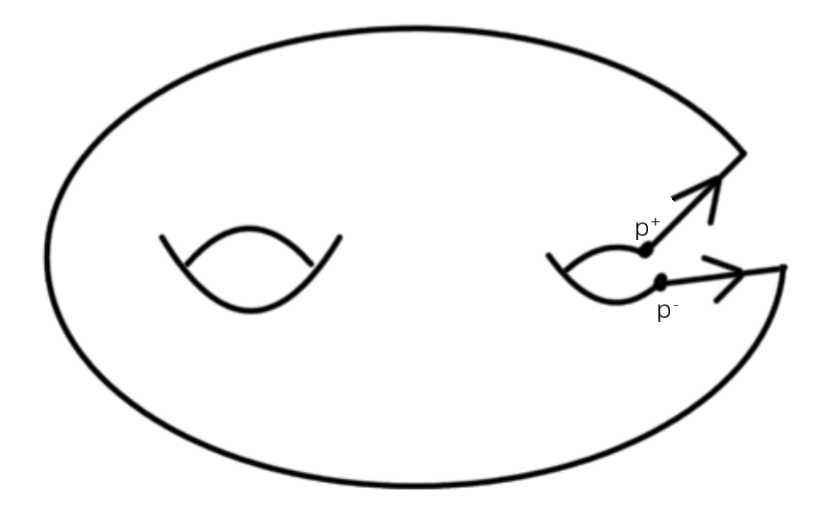}
    \caption{Pinching}
    \label{pinching}
\end{figure}

%Now we get $F_3$
%So we get , which is from some full measure subset of $ \H_g(2^{g-1})$ to $\H_{g-1}(2^{g-1},0,0)\times A_{[\frac{a}{\sqrt{g}},\frac{b}{\sqrt{g}}]}$.
\subsection{Collapsing a closed saddle connection}\label{F_3}
Now denote by
$$\hat{\mathcal{C}}^{1,0}_{g,Per}(2^{g-1})[\frac{a}{\sqrt{g}},\frac{b}{\sqrt{g}}]$$
the set of $(X,\omega,\gamma)$, where $(X,\omega)\in \mathcal{C}^{1,0}_{g,Per}(2^{g-1})[\frac{a}{\sqrt{g}},\frac{b}{\sqrt{g}}]$ and $\gamma$ is
a closed saddle connection on $(X,\omega)$ with length in $[\frac{a}{\sqrt{g}},\frac{b}{\sqrt{g}}]$. 
Then for $(X,\omega,\gamma)\in \hat{\mathcal{C}}^{1,0}_{g,Per}(2^{g-1})[\frac{a}{\sqrt{g}},\frac{b}{\sqrt{g}}]$, we make opening, moving and pinching surgeries in sequence
and get a new surface in $\H_{g-1}(2^{g-2},0,0)$ and a vector in $A_{[\frac{a}{\sqrt{g}},\frac{b}{\sqrt{g}}]}$. 
This defines a mapping
$$F:\hat{\mathcal{C}}^{1,0}_{g,Per}(2^{g-1})[\frac{a}{\sqrt{g}},\frac{b}{\sqrt{g}}] \to |M_1|\H_{g-1}(2^{g-2},0,0)\times A_{[\frac{a}{\sqrt{g}},\frac{b}{\sqrt{g}}]},$$
where $|M_1|=g-1$ denote which zero we have collapsed.
Next we will show the mapping is one to one and measure-preserving almost everywhere. 

\begin{proposition}
$F$ is invertible except for a subset of measure zero of 
$$|M_1|\H_{g-1}(2^{g-2},0,0)\times A_{[\frac{a}{\sqrt{g}},\frac{b}{\sqrt{g}}]}.$$
\end{proposition}
\begin{proof}
Let $[(X,\omega),\vec{a}]\in \H_{g-1}(2^{g-2},0,0)\times A_{[\frac{a}{\sqrt{g}},\frac{b}{\sqrt{g}}]}$.
Since the opening surgery can be inversed, it suffices to consider when the pinching and moving surgeries can be inversed.

For pinching surgery, it can be inversed except for some directions which have saddle connections, which is a measure zero subset $A_0(X,\omega)$
in $A_{[\frac{a}{\sqrt{g}},\frac{b}{\sqrt{g}}]}$ for every surface in $\H_{g-1}(2^{g-2},0,0)$: 
choose a vector $\kappa$ in $A_{[\frac{a}{\sqrt{g}},\frac{b}{\sqrt{g}}]}$, 
consider the two rays from $p^+$ and $p^-$ whose holonomy is $\kappa$, if the rays exist, then cut the surface along the two rays
to obtain two closed loci and glue them. 

Suppose after the inverse of pinching $[(X,\omega),\vec{a}]$ becomes $(X',\omega',\gamma)$, 
where $\gamma$ is a simple closed curve consisting of two saddle connections with same length.
For moving surgery, from \ref{F_2}, it can be inversed except for $(X',\omega')$ having the other closed curve $\gamma'$ consisting of a pair of non-homologous 
saddle connections which have same length. But $\gamma'$ is fixed during the pinching surgery and other saddle connections are not changed,
this means $(X,\omega)$ having a simple closed curve consisting of two saddle connections with same length,
which is a measure zero set of $\H_{g-1}(2^{g-2},0,0)$.

So for any $[(X,\omega),\vec{a}]\in \H_{g-1}(2^{g-2},0,0)\times A_{[\frac{a}{\sqrt{g}},\frac{b}{\sqrt{g}}]}$ except for a subset of measure zero, the surgeries can be inversed,
note that after the inversed surgery the new zero we get has $g-1$ choices to be marked,
which means the surgery is $(g-1)$-to-$1$ and $F$ is invertible almost everywhere on 
$|M_1|\H_{g-1}(2^{g-2},0,0)\times A_{[\frac{a}{\sqrt{g}},\frac{b}{\sqrt{g}}]}.$
\end{proof}

Now we will see the mapping is also measure-preserving on $\mathcal{C}^{1,0}_{g,Per}(2^{g-1})[\frac{a}{\sqrt{g}},\frac{b}{\sqrt{g}}]$.
In fact, under local chart determined by the basis we choose,
$$F(z_1,\cdots,z_n)=(z_1,z_i+\upsilon(z_i)).$$ 
Here $z_1$ means the holonomy of the saddle connection $\gamma$, which becomes a
vector in the annulus.
The Jacobian of $F$ has determinant one. Since Masur-Veech measure is invariant under coordinate transformation of period mapping, $F$ is measure-preserving locally.
Let $\tilde{\H}$ be the subset on which $F $ is invertible,
which is a full measure subset of $\H_{g-1}(2^{g-2},0,0)\times A_{[\frac{a}{\sqrt{g}},\frac{b}{\sqrt{g}}]}$.
Then we have
$$|M_1|\tilde{\H} \subset F( \hat{\mathcal{C}}^{1,0}_{g,Per}(2^{g-1}))\subset |M_1|\H_{g-1}(2^{g-2},0,0)\times A_{[\frac{a}{\sqrt{g}},\frac{b}{\sqrt{g}}]}.$$

Since $F$ is measure-preserving locally and $\H_{g-1}(2^{g-2},0,0)\times A_{[\frac{a}{\sqrt{g}},\frac{b}{\sqrt{g}}]}\setminus \tilde{\H}$ is a measure-zero set,
$F^{-1}(\H_{g-1}(2^{g-2},0,0)\times A_{[\frac{a}{\sqrt{g}},\frac{b}{\sqrt{g}}]}\setminus \tilde{\H})$ is also a measure-zero set.
So we have
$$\V(\hat{\mathcal{C}}^{1,0}_{g,Per}(2^{g-1})[\frac{a}{\sqrt{g}},\frac{b}{\sqrt{g}}])=(g-1)\V(F^{-1}(\tilde{\H}))$$

Moreover, $F$ is measure-preserving on $F^{-1}(\tilde{\H})$, so we have
$$\V(F^{-1}(\tilde{\H}))=\V(\tilde{\H})=\V([\H_{g-1}(2^{g-2},0,0)\times (A_{[\frac{a}{\sqrt{g}},\frac{b}{\sqrt{g}}]})])$$.
Combine the two equalities above, 
$$\V(\hat{\mathcal{C}}^{1,0}_{g,Per}(2^{g-1})[\frac{a}{\sqrt{g}},\frac{b}{\sqrt{g}}])=(g-1)\V([\H_{g-1}(2^{g-2},0,0)\times (A_{[\frac{a}{\sqrt{g}},\frac{b}{\sqrt{g}}]})]).$$

%$$F_3( \mathcal{C}^{1,0}_{g,Per}(2^{g-1}))=M\times [\H_{g-1}(2^{g-2},0,0)\times (A_{[\frac{a}{\sqrt{g}},\frac{b}{\sqrt{g}}]})]\setminus \tilde{\H},$$
%Since $\V(\H_{g-1}(2^{g-2},0,0))=\V(\H_{g-1}(2^{g-2}))$, $|M|=g-1$, let $g\to \infty$ and by \ref{Volume} we have

%$$\V(\mathcal{C}^{1,0}_{g}(2^{g-1})[\frac{a}{\sqrt{g}},\frac{b}{\sqrt{g}}])=\V(\mathcal{C}^{1,0}_{g,Per}(2^{g-1})[\frac{a}{\sqrt{g}},\frac{b}{\sqrt{g}}]).$$

%\begin{proposition}
%$$\lim_{g\to \infty}\frac{\mathcal{C}^{1,0}_{g,Per}(2^{g-1})}{\V(\H_{g}(2^{g-1}))}=3\pi(b^2-a^2).$$
%\end{proposition}

\section{Collapse closed saddle connections simultaneously}
In this section we will prove Theorem \ref{closed}.

From the construction above we can collapse one closed saddle connection on $(X,\omega)$ in a nearly full measure set of $\mathcal{C}^{1,0}_{g}(2^{g-1})$ and obtain a new surface
$(X',\omega')\in \H_{g-1}(2^{g-2},0,0)$. If we want to collapse any $k$ closed saddle connections, 
we have to ensure the locus of surgery for these saddle connections are disjoint.
For a given interval $[\frac{a}{\sqrt{g}},\frac{b}{\sqrt{g}}]$,
define the \emph{exception set} $\mathcal{C}^{exc}_{g,Per}(2^{g-1})[\frac{a}{\sqrt{g}},\frac{b}{\sqrt{g}}]$ to be the subset of 
$\mathcal{C}^{1,0}_{g,Per}(2^{g-1})[\frac{a}{\sqrt{g}},\frac{b}{\sqrt{g}}]$ 
such that on every $(X,\omega)\in \mathcal{C}^{exc}_{g,Per}(2^{g-1})[\frac{a}{\sqrt{g}},\frac{b}{\sqrt{g}}]$ there are two closed saddle connections with length in
$[\frac{a}{\sqrt{g}},\frac{b}{\sqrt{g}}]$ whose surgery loci intersect.
First we will show that the measure of exception set goes to zero as genus goes to infinity.

\subsection{The measure of exception set}\label{exception}

For $(X,\omega)\in \mathcal{C}^{exc}_{g,Per}(2^{g-1})[\frac{a}{\sqrt{g}},\frac{b}{\sqrt{g}}]$, there are two situations: 
\begin{itemize}
  \item[1]  There exist two closed saddle connections on $(X,\omega)$ with lengths in $[\frac{a}{\sqrt{g}},\frac{b}{\sqrt{g}}]$
            that do not share a zero. And the loci of opening surgery intersect or the two closed saddle connections intersect.
  \item[2]  There exist two closed saddle connections on $(X,\omega)$ with lengths in $[\frac{a}{\sqrt{g}},\frac{b}{\sqrt{g}}]$ sharing one zero.
\end{itemize}

Note that in the first situation, we can find a curve connecting the two zeros with length no more than $\frac{2b}{\sqrt{g}}$.

Fix $B\in \mathbb{R}^+$, define $\mathcal{C}_g^{1}(\frac{B}{\sqrt{g}})$ to be the set of $(X,\omega)\in \mathcal{C}^{1,0}_{g,Per}(2^{g-1})[0,\frac{B}{\sqrt{g}}]$ 
such that there exist two closed saddle connections on $(X,\omega)$
that do not share a zero, but the lengths of the two closed saddle connections and the distance between the two zeros are less than $\frac{B}{\sqrt{g}}$.
Define $\mathcal{C}_g^{2}(\frac{B}{\sqrt{g}})$ to be the set of $(X,\omega)\in \mathcal{C}^{1,0}_{g,Per}(2^{g-1})[0,\frac{B}{\sqrt{g}}]$ 
such that there exist two closed saddle connections with lengths less than $\frac{B}{\sqrt{g}}$ sharing one zero.
Obviously for $B\geq 2b$, we have 
$$ \mathcal{C}^{exc}_{g,Per}(2^{g-1})[\frac{a}{\sqrt{g}},\frac{b}{\sqrt{g}}]\subset \mathcal{C}_g^{1}(\frac{B}{\sqrt{g}}) \cup \mathcal{C}_g^{2}(\frac{B}{\sqrt{g}}).$$

We need to compute the measure of $\mathcal{C}_g^{1}(\frac{B}{\sqrt{g}})$ and $\mathcal{C}_g^{2}(\frac{B}{\sqrt{g}})$ when $g\to \infty$.
First consider $\mathcal{C}_g^{2}(\frac{B}{\sqrt{g}})$, we have
\begin{proposition}\label{chain}
$$\lim_{g\to \infty}\frac{\V(\mathcal{C}_g^{2}(\frac{B}{\sqrt{g}}))}{\V(\H_g(2^{g-1}))} = 0.$$
\end{proposition}

\begin{proof}
For $(X,\omega) \in \mathcal{C}_g^{2}(\frac{B}{\sqrt{g}})$, 
let $p$ be a zero on $(X,\omega)$ and $\gamma_1, \gamma_2$ be the closed saddle connections with lengths less than $\frac{B}{\sqrt{g}}$ at $p$.

As above we can collapse $\gamma_1$ and get a new translation surface in $\H_{g-1}(2^{g-2},0,0)$. 
After the surgery $\gamma_2$ will become a segment connecting two marked regular points which are from the pinching surgery.
So the image of mapping $F_3$ on $\mathcal{C}_g^{2}(\frac{B}{\sqrt{g}})$ is in 
$\H_{g-1,\frac{B}{\sqrt{g}}}(2^{g-2},0,0)\times D(\frac{B}{\sqrt{g}})$, 
where $\H_{g-1,\frac{B}{\sqrt{g}}}(2^{g-2},0,0)$ is the subset that there exists a segment connecting the two marked points of length less than 
$\frac{B}{\sqrt{g}}$ and $D(\frac{B}{\sqrt{g}})$ is the disk of radius $\frac{B}{\sqrt{g}}$.
From \cite[Theorem 1.2]{aggarwal2019large}, the Siegel-Veech constant of saddle connections connecting two fixed zeros of order $m_1$ and $m_2$ is
$$c=(m_1+1)(m_2+1)(1+O(\frac{1}{g})).$$
So we have
$$\V(\mathcal{C}_g^{2}(\frac{B}{\sqrt{g}}))\leq (g-1)\frac{B^2}{g}\V(\H_{g-1,\frac{B}{\sqrt{g}}}(2^{g-2},0,0)) \leq c(g-1)\frac{B^4}{g^2}\V(\H_{g-1}(2^{g-2},0,0)),$$
where the coefficient $g-1$ is the choice of $p$.
Then by \ref{Volume}, we have
$$\lim_{g\to \infty}\frac{\V(\mathcal{C}_g^{2}(\frac{B}{\sqrt{g}}))}{\V(\H_g(2^{g-1}))}\leq \lim_{g\to \infty}3c\frac{B^4}{g}=0$$
\end{proof}

For $\mathcal{C}_g^{1}(\frac{B}{\sqrt{g}})$ we also have
\begin{proposition}\label{intersect}
$$\lim_{g\to \infty}\frac{\V(\mathcal{C}_g^{1}(\frac{B}{\sqrt{g}}))}{\V(\H_g(2^{g-1}))} = 0.$$
\end{proposition}

\begin{proof}
For $(X,\omega) \in \mathcal{C}_g^{1}(\frac{B}{\sqrt{g}})$, let $p_1, p_2$ be the two zeros on $(X,\omega)$ with distance less than $\frac{B}{\sqrt{g}}$ 
and $\gamma_1, \gamma_2$ be the two closed saddle connections connecting $p_1, p_2$.
For every $k\in \mathbb{N}^+$, denote by $\mathcal{C}_g^{1,k}(\frac{B}{\sqrt{g}})$ the subset of $\mathcal{C}_g^{1}(\frac{B}{\sqrt{g}})$
that the geodesic between $p_1$ and $p_2$ is a concatenation of $k$ open saddle connections.
Their lengths are less than $\frac{B}{\sqrt{g}}$.
Choose the shortest saddle connection, using the surgery in \cite{masur2024lengths} we can collapse it and get a new translation surface in 
$\H_g(3,1,2^{g-3})$. This is because the locus, that is the extension of the shortest saddle connection cannot intersect the whole geodesic.
Repeat the surgery until $p_1$ and $p_2$ are collapsed to one zero and we get a map
$$F:\mathcal{C}_g^{1,k}(\frac{B}{\sqrt{g}})\to \mathcal{C}_g^{2}(2+k,1^k,2^{g-2-k})(\frac{B}{\sqrt{g}})\times D^k(\frac{B}{\sqrt{g}}),$$
where $\mathcal{C}_g^{2}(2+k,1^k,2^{g-2-k})(\frac{B}{\sqrt{g}})$ is the subset of $\H_g(2+k,1^k,2^{g-2-k})$ such that there are
two closed saddle connections $\gamma_1, \gamma_2$ connecting the zero of order $2+k$ with lengths less than $\frac{B}{\sqrt{g}}$ .

Similarly one can collapse $\gamma_1$ and the other will become a saddle connection connecting two zeros. 
Note that the order is larger than $2$, although the collapsing surgery can also be realized, we need to require the angles 
on the two sides of the closed saddle connections collapsed, which decide the order of zeros after collapsing. 
We will explain this situation in detail in Section \ref{general case}.

Let $\mathcal{C}_g^{2,b',b''}(2+k,1^k,2^{g-2-k})(\frac{B}{\sqrt{g}})$ be the subset of 
$\mathcal{C}_g^{2}(2+k,1^k,2^{g-2-k})(\frac{B}{\sqrt{g}})$ such that 
the two angles on both sides of $\gamma_1$ are $(2b'+1)\pi$ and $(2b''+1)\pi$, where $b'+b''=2+k$.
Then after collapsing we get a map
$$F':\mathcal{C}_g^{2,b',b''}(2+k,1^k,2^{g-2-k})(\frac{B}{\sqrt{g}}) \to \H_{g-1,\frac{B}{\sqrt{g}}}(b'-1,b''-1,1^k,2^{g-2-k})\times D(\frac{B}{\sqrt{g}}).$$

Similar to Proposition \ref{chain}, we have
$$\V(\mathcal{C}_g^{2,b',b''}(2+k,1^k,2^{g-2-k})(\frac{B}{\sqrt{g}}))\leq b'b''\frac{B^4}{g^2}\V(\H_{g-1}(b'-1,b''-1,1^k,2^{g-2-k})).$$
Summing over all $k$ and $(b',b'')$ we have
$$\V(\mathcal{C}_g^{1}(\frac{B}{\sqrt{g}}))\leq (g-1)(g-2)\sum_k \sum_{(b',b'')}b'b''\frac{B^4}{g^2}\V(\H_{g-1}(b'-1,b''-1,1^k,2^{g-2-k}))(\frac{B^2}{g})^k,$$
where the coefficient $(g-1)(g-2)$ is the choice of $p_1$ and $p_2$.
Then when $g\to \infty$ we have
$$\frac{\V(\mathcal{C}_g^{1}(\frac{B}{\sqrt{g}}))}{\V(\H_g(2^{g-1}))}\leq
\sum_k \sum_{(b',b'')}b'b''B^4\frac{\V(\H_{g-1}(b'-1,b''-1,1^k,2^{g-2-k}))}{\V(\H_g(2^{g-1}))}(\frac{B^2}{g})^k=O(\frac{1}{g}).$$
\end{proof}

Combine the two propositions we have
\begin{corollary}
$$\frac{\V(\mathcal{C}^{exc}_{g,Per}(2^{g-1})[\frac{a}{\sqrt{g}},\frac{b}{\sqrt{g}}])}{\V(\H_g(2^{g-1}))}=O(\frac{1}{g}).$$
\end{corollary}

\subsection{Collapsing simultaneously}
Now fix a positive integer $K$ and $n_1,\cdots,n_k$ a partition of $K$.
We want to compute the limit of the expectation
$$(L_{g,[a_1,b_1]})_{n_1}\cdots (L_{g,[a_k,b_k]})_{n_k},$$
for which we will collapse $K$ closed saddle connections in order.

Define $\hat{\H_g}(2^{g-1})=(X,\omega,\Gamma_1,\cdots,\Gamma_k)$, 
where $(X,\omega)\in \H_g(2^{g-1})$, 
and $\Gamma_i$ is an ordered list of $n_i$ oriented closed saddle connections $(\gamma_{1,i},\cdots,\gamma_{n_i,i})$ on $(X,\omega)$
with lengths in $[\frac{a_i}{\sqrt{g}},\frac{b_i}{\sqrt{g}}]$. 
Suppose $a_1\leq b_1 \leq a_2\leq b_2 \cdots \leq a_k \leq b_k$, let $B=b_k$.
Then we have
$$\bigcup_i \mathcal{C}^{exc}_{g,Per}(2^{g-1})[\frac{a_i}{\sqrt{g}},\frac{b_i}{\sqrt{g}}]\subset \mathcal{C}_g^{1}(\frac{B}{\sqrt{g}})\cup \mathcal{C}_g^{2}(\frac{B}{\sqrt{g}}).$$

Let $\H'_g(2^{g-1})=\H_g(2^{g-1})\setminus (\mathcal{C}_g^{1}(\frac{B}{\sqrt{g}})\cup \mathcal{C}_g^{2}(\frac{B}{\sqrt{g}}))$, and define
$\hat{\H}'_g(2^{g-1})$ as $\hat{\H_g}(2^{g-1})$.
We will adjust the surgery from last section on permissible set to collapse the closed saddle connections in
$\Gamma_1,\cdots,\Gamma_k$ in order.
The opening and pinching operations are same to the situation of one curve. The only modification is the moving surgery.

\begin{remark}
Note that if $K=1$, $\hat{\H}'_g(2^{g-1})$ is just the set $\hat{\mathcal{C}}^{1,0}_{g,Per}(2^{g-1})[\frac{a}{\sqrt{g}},\frac{b}{\sqrt{g}}]$ in Section \ref{F_3}.
\end{remark}

\iffalse
Since we can only consider the curves that go through
at most $g^{2/3}$ handles, all curves in $\Gamma_1,\cdots,\Gamma_k$ 
stay in a subsurface of at most $Kg^{2/3}$.
We choose the first $2K$ basis to be $\alpha_1,\cdots\alpha_K$ and their dual $\beta_1\cdots \beta_K$. 
For $\beta_{K+1}$, consider the set consisting of $f\circ h^{m_i}\circ f^{-1}(\alpha_j)$
that are disjoint with $\alpha_1,\cdots\alpha_K$.
Choose $\beta_{K+1}$ in the set whose holonomy has minimal vertical projection,
and choose its dual curve to be $\alpha_{K+1}$.
Then we choose the tangent vector $\upsilon$ satisfying $\upsilon(\gamma_i)=-\omega(\gamma_i)$, and adjust the value of $\upsilon(\alpha_{2K+1})$ to ensure the area unchanged.
Since $K$ is fixed and $g\to \infty$, this $\upsilon$ can be chosen.
\fi
For $K$ curves, consider its preimages under the marking, we can also use the Euclid's Extended Algorithm to find their dual curves. Then consider the symplectic complement
of the subspace generated by these curves and repeat the process to extend a symplectic basis on the surface. Finally we can define the moving vector $\upsilon$ and carry out 
the moving surgery.

Finally  we obtain a mapping
$$\hat{F}: \hat{\H}'_g(2^{g-1}) \to M_K\cdot \H_{g-K}(2^{g-1-K},0^{2K}) \times \prod^k_{i=1} (A_{[\frac{a_i}{\sqrt{g}},\frac{b_i}{\sqrt{g}}]})^{n_i},$$
where $\hat{\H}'_g(2^{g-1})$ is defined as $\hat{\H}_g(2^{g-1})$
and $M_K$ is a combinatorial datum consisting of choosing $K$ ordered zeros from $g-1$ zeros to label. 
From section \ref{F_3}, this map is one to one on a total measure subset of $\H_{g-K}(2^{g-1-K},0^{2K})$, 
since at a marked point there is only one choice to split along a given vector.

As the discussion above, 
we have 
%$$\V(\hat{\mathcal{C}}^{1,0}_{g,Per}(2^{g-1})\setminus (\hat{\mathcal{C}}_g^{1}(\frac{B}{\sqrt{g}})\cup \hat{\mathcal{C}}_g^{2}(\frac{B}{\sqrt{g}})))=|M_K|$$
\begin{equation*}
\begin{aligned}
&\int_{\hat{\H}'_g(2^{g-1})}(L_{g,[a_1,b_1]})_{n_1}\cdots (L_{g,[a_k,b_k]})_{n_k}d\mu_{MV}\\&=\frac{1}{2^K}\V(\hat{\H}'_g(2^{g-1}))
\\&=|M_K|\prod^k_{i=1}[\frac{\pi(b^2-a^2)}{2g}]^{n_i}\V(\H_{g-k}(2^{g-1-K}))\\
%&\to g^K\prod^k_{i=1}[\frac{\pi(b_i^2-a_i^2)}{g}]^{n_i}\V(\H_{g-K}(2^{g-1}))\\
&\to \prod^k_{i=1}[\frac{\pi}{2}(b_i^2-a_i^2)]^{n_i}\V(\H_{g-K}(2^{g-1-K})), g\to \infty,
\end{aligned}
\end{equation*}
where the factor $\frac{1}{2^K}$ is because we omit the orientation when counting the saddle connections
and the limit is because $|M_K|=(g-1)\cdots (g-K)$ and
$$\lim_{g\to \infty}\frac{(g-1)\cdots (g-K)}{g^{n_1+\cdots+n_k}}=\lim_{g\to \infty}\frac{(g-1)\cdots (g-K)}{g^K}=1.$$

Let $L'_{g,[a_i,b_i]}$ be the restriction of $L_{g,[a_i,b_i]}$ on $\hat{\H}'_g(2^{g-1})$.
Then the limit of factorial moment of $(L'_{g,[a_i,b_i]})_i$ is
$$\lim_{g\to \infty}\mathbb{E}[(L'_{g,[a_1,b_1]})_{n_1}\cdots (L'_{g,[a_k,b_k]})_{n_k}]=\lim_{g\to \infty}\frac{\V(\H_{g-K}(2^{g-1-K}))}{2^K\V(\hat{\H}'_g(2^{g-1}))}\prod^k_{i=1}[\pi(b_i^2-a_i^2)]^{n_i}.$$
From Proposition \ref{chain} and \ref{intersect}, we have
$$\lim_{g\to \infty}\frac{\V(\hat{\H}'_g(2^{g-1}))}{\V(\H_g(2^{g-1}))}=1.$$
From Theorem \ref{Volume}, we have
$$\lim_{g\to \infty}\frac{\V(\H_{g-K}(2^{g-1-K}))}{\V(\H_{g}(2^{g-1}))}=3^K.$$
So we have
$$\lim_{g\to \infty}\mathbb{E}[(L'_{g,[a_1,b_1]})_{n_1}\cdots (L'_{g,[a_k,b_k]})_{n_k}]=\prod^k_{i=1}[\frac{3\pi}{2}(b_i^2-a_i^2)]^{n_i}.$$
By theorem \ref{moment}, $(L'_{g,[a_1,b_1]},\cdots,L'_{g,[a_k,b_k]})$ converges jointly in distribution to a vector of random variables with Poisson distributions of means $\lambda_{[a_i,b_i]}$, where
$$\lambda_{[a_i,b_i]}=\frac{3\pi}{2}(b_i^2-a_i^2).$$

From the property of Poisson distribution, we have
$$\mathbb{P}(L'_{g,[a_1,b_1]}=n_1,\cdots,L'_{g,[a_k,b_k]}=n_k)=\prod_{i=1}^k\frac{\lambda_{[a_i,b_i]}^{n_i}e^{-\lambda_{[a_i,b_i]}}}{n_i!}.$$
On the other hand 
$$\mathbb{P}(L'_{g,[a_1,b_1]}=n_1,\cdots,L'_{g,[a_k,b_k]}=n_k)=\frac{\V(\{(X,\omega):L'_{g,[a_i,b_i]}(X,\omega)=n_i,i=1,\cdots,k\})}{\V(\hat{\H}'_g(2^{g-1}))}.$$
Note that 
\begin{equation*}
\begin{aligned}
&\V(\{(X,\omega):L_{g,[a_i,b_i]}(X,\omega)=n_i,i=1,\cdots,k\})-
\V(\mathcal{C}_g^{1}(\frac{B}{\sqrt{g}})\cup \mathcal{C}_g^{2}(\frac{B}{\sqrt{g}}))\\ &\leq \V(\{(X,\omega):L'_{g,[a_i,b_i]}(X,\omega)=n_i,i=1,\cdots,k\})\\&\leq
\V(\{(X,\omega):L_{g,[a_i,b_i]}(X,\omega)=n_i,i=1,\cdots,k\}).
\end{aligned}
\end{equation*}
Since
$$\lim_{g\to \infty}\frac{\V(\hat{\H}'_g(2^{g-1}))}{\V(\H_g(2^{g-1}))}=1,\\ \lim_{g\to \infty}\frac{\V(\mathcal{C}_g^{1}(\frac{B}{\sqrt{g}})\cup \mathcal{C}_g^{2}(\frac{B}{\sqrt{g}}))}{\V(\H_g(2^{g-1}))}=0.$$
We have
\begin{equation*}
\begin{aligned}
&\mathbb{P}(L_{g,[a_1,b_1]}=n_1,\cdots,L_{g,[a_k,b_k]}=n_k)-
\frac{\V(\mathcal{C}_g^{1}(\frac{B}{\sqrt{g}})\cup \mathcal{C}_g^{2}(\frac{B}{\sqrt{g}}))}{\V(\H_g(2^{g-1}))}\\ &\leq \mathbb{P}(L'_{g,[a_1,b_1]}=n_1,\cdots,L'_{g,[a_k,b_k]}=n_k)\frac{\V(\hat{\H}'_g(2^{g-1}))}{\V(\H_g(2^{g-1}))}\\&\leq
\mathbb{P}(L_{g,[a_1,b_1]}=n_1,\cdots,L_{g,[a_k,b_k]}=n_k).
\end{aligned}
\end{equation*}
Let $g \to \infty$ we have
$$\mathbb{P}(L_{g,[a_1,b_1]}=n_1,\cdots,L_{g,[a_k,b_k]}=n_k)=\mathbb{P}(L'_{g,[a_1,b_1]}=n_1,\cdots,L'_{g,[a_k,b_k]}=n_k).$$
So $(L_{g,[a_1,b_1]},\cdots,L_{g,[a_k,b_k]})$ also converges jointly in distribution to a vector of random variables 
with Poisson distributions of means $\lambda_{[a_i,b_i]}$, where
$$\lambda_{[a_i,b_i]}=\frac{3\pi}{2}(b_i^2-a_i^2).$$

Then by Theorem \ref{moment}, Theorem \ref{closed} is proved.
\section{General stratum}\label{general case}
In this section, we first consider the stratum $\H_g(m^{O(g)},1^{2g-2-mO(g)})$, where $m\geq 3$ and $\lim\limits_{g\to \infty}\frac{O(g)}{g}=c$.
If a stratum has zeros of order more than $3$, the angles at both sides are not certain. 
But we can consider some fixed configuration of closed saddle connection.
For a configuration $(J,b'_k,b''_k,a'_i,a''_i)$ of closed saddle connection,
We know it suffices to consider multiplicity $1$ without cylinders, then the configuration becomes $(b',b'')$, where $b'+b''=m$, 
which means the angles of the closed saddle connection are $(2b'+1)\pi$ on one side and $(2b''+1)\pi$ on the other side.
Now we choose all such configurations at each zeros of order $m$.

For $k$ such closed saddle connections, we can use the surgery in Section \ref{surgery} to collapse them. 
Note that in this situation the opening surgery will replace the zero of order $m$ by two zeros of order $1$ and $m-1$, 
and the angles at both sides of the zero of $m-1$ order are $(2b')\pi$ and $(2b'')\pi$. Then after moving and pinching, 
the simple zero becomes two regular points and the zero of $m-1$ order becomes two zeros of order $b'-1$ and $b''-1$.
Finally we obtain a mapping
$$\hat{F}: \hat{\mathcal{H}}'_{g}(m^{O(g)},1^{O'(g)})\to M_K\cdot \H_{g-K}(\kappa')\times \prod^k_{i=1} (A_{[\frac{a_i}{\sqrt{g}},\frac{b_i}{\sqrt{g}}]})^{n_i},$$
where $\kappa'=((b'-1)^K,(b''-1)^K,m^{O(g)-K},1^{O'(g)})$.

In this situation, the construction of the inverse mapping is a little different.
Suppose we have chosen a fixed vector sequence $(v_1,\cdots,v_K)\in \prod^k_{i=1} (A_{[\frac{a_i}{\sqrt{g}},\frac{b_i}{\sqrt{g}}]})^{n_i}$, 
for some $v_i$, we want to split the surface in $\H_{g-K}((b'-1)^K,(b''-1)^K,m^{O(g)-K},1^{O'(g)})$ along $v_i$ at a pair of zeros of orders $(b'-1)$ and $(b''-1)$,
but we have $b'b''$ choices for the split locus. So the mapping $\hat{F}$ is $(b'b'')^K$ to one.
So we have

\begin{equation*}
\begin{aligned}
\V(\hat{\mathcal{H}}'_{g}(m^{O(g)},1^{O'(g)}))=|M_K|(b'b'')^K \V(\H_{g-K}(\kappa'))\prod^k_{i=1}(\pi \frac{b_i^2-a_i^2}{g})^{n_i},
\end{aligned}
\end{equation*}

where $|M_K|=O(g)(O(g)-1)\cdots(O(g)-K+1)$.
And similar to the proof of theorem \ref{closed} we have
\begin{equation*}
\begin{aligned}
&\lim_{g\to \infty}\mathbb{E}[(L'_{g,[a_1,b_1]})_{n_1}\cdots (L'_{g,[a_k,b_k]})_{n_k}]
\\&=|M_K|(b'b'')^K\frac{\V(\H_{g-K}(\kappa'))}{2^K\V(\H_g(m^{O(g)},1^{O'(g)}))}\prod^k_{i=1}(\pi \frac{b_i^2-a_i^2}{g})^{n_i}
\\&=(b'b'')^K \frac{(m+1)^K}{2^K(b'b'')^K} \prod^k_{i=1}[\pi(b_i^2-a_i^2)]^{n_i}\frac{O(g)(O(g)-1)\cdots(O(g)-K)}{g^K}
\\&=c^K(\frac{m+1}{2})^K\prod^k_{i=1}[\pi(b_i^2-a_i^2)]^{n_i}
\\&=\prod^k_{i=1}[\frac{c(m+1)\pi}{2}(b_i^2-a_i^2)]^{n_i}.
\end{aligned}
\end{equation*}

So we have
\begin{theorem}
Suppose $\lim\limits_{g\to \infty}\frac{O(g)}{g}=c$. 
For the stratum $\H(m^{O(g)},1^{2g-2-mO(g)})$, a configuration $\mathcal{C}=(b',b'')$ and $a,b\in \mathbb{R}$,
let the random variable
$$L_{\mathcal{C},g,[a,b]}:\H(m^{O(g)},1^{2g-2-mO(g)}) \to \mathbb{N}_0$$ 
be the number of closed saddle connections on $(X,\omega)$ satisfying the configuration with lengths in $[\frac{a}{\sqrt{g}},\frac{b}{\sqrt{g}}]$.
Then for disjoint intervals $[a_1,b_1],\cdots,[a_k,b_k]$, when $g\to \infty$, the random variable sequence
$$(L_{\mathcal{C},g,[a_1,b_1]},\cdots,L_{\mathcal{C},g,[a_k,b_k]}):\H(m^{O(g)},1^{2g-2-mO(g)}) \to \mathbb{N}^k_0$$
converges jointly in distribution to a vector of random variables with Poisson distributions of means $\lambda_{[a_i,b_i]}$, where
$$\lambda_{[a_i,b_i]}=\frac{c}{2}(m+1)\pi(b_i^2-a_i^2).$$

\end{theorem}

Note that $\lambda_{[a_i,b_i]}$ is independent of $(b',b'')$,
now we consider general closed saddle connections in $\H(m^{O(g)},1^{2g-2-mO(g)})$.
Let 
$$L_{g,[a,b]}:\H(m^{O(g)},1^{2g-2-mO(g)}) \to \mathbb{N}_0$$
be the number of closed saddle connections on $(X,\omega)$ with lengths in $[\frac{a}{\sqrt{g}},\frac{b}{\sqrt{g}}]$.
Then we have 
$$L_{g,[a,b]}=\sum_{i=1}^{m-1}L_{\mathcal{C}_i,g,[a,b]},$$
where $\mathcal{C}_i=(i,m-i)$. 

Fix some $i,j\in \{1,2,\cdots,m-1\}$,
let
$$\hat{L}^{i,j}_{g,[a,b]}=(L_{\mathcal{C}_i,g,[a,b]},L_{\mathcal{C}_j,g,[a,b]}).$$
We will compute the distribution of $\hat{L}^{i,j}_{g,[a,b]}$.
For $n_i,n_j$, consider the set consisting of 
$$(X,\omega,\Gamma_i,\Gamma_j),$$
where $(X,\omega)\in \mathcal{C}^{1,0}_{g}(m^{O(g)},1^{2g-2-mO(g)})$, 
and $\Gamma_i$ is an ordered list of $n_i$ closed saddle connections $(\gamma_{1,i},\cdots,\gamma_{n_i,i})$ 
with lengths in $[\frac{a}{\sqrt{g}},\frac{b}{\sqrt{g}}]$ and configuration $\mathcal{C}_i$.

If we want to collapse curves in $\Gamma_i$ and $\Gamma_j$, we need to consider the exception set where the surgery loci intersect.
Similar to \ref{exception}, if two surgery loci intersect, then either the two zeros are close or the two closed saddle connections share
the same zero.
From section \ref{exception}, the measure of such subset goes to zero as $g\to \infty$.
So we can also collapse the curves in $\Gamma_i$ and $\Gamma_j$ simultaneously.
Again we have
\begin{proposition}
$\hat{L}^{i,j}_{g,[a,b]}$ converges jointly in distribution to a vector of random variables with Poisson distributions of means $\lambda_{[a,b]}$, where
$$\lambda_{[a,b]}=\frac{c}{2}(m+1)\pi(b^2-a^2).$$

This means for different configurations, the limits of $L_{\mathcal{C}_i,g,[a,b]}$ are independent.
\end{proposition}

From the property of Poisson distribution, we have
\begin{proposition}\label{order m}
$L_{g,[a,b]}=\sum\limits_{i=1}^{m-1}L_{\mathcal{C}_i,g,[a,b]}$ converges in distribution to a Poisson distribution of mean $\lambda_{[a,b]}$, where
$$\lambda_{[a,b]}=\frac{c}{2}(m^2-1)\pi(b^2-a^2).$$
\end{proposition}

\begin{proof}
Since $L_{\mathcal{C}_i,g,[a,b]}$ are independent and 
$$L_{g,[a,b]}=\sum_{i=1}^{m-1}L_{\mathcal{C}_i,g,[a,b]}.$$

We have 
$L_{g,[a,b]}$ converges in distribution to a Poisson distribution, and the mean is
$(m-1)\lambda_{[a,b]}=c(m^2-1)\pi(b^2-a^2)$.
\end{proof}

Finally, consider $\H(m_1^{O_1(g)},\cdots,m_k^{O_k(g)},1^{O(g)})$, where $\lim\limits_{g\to \infty}\frac{O_i(g)}{g}=c_i.$
Let 
$$L_{g,i,[a,b]}:\H(m_1^{O_1(g)},\cdots,m_k^{O_k(g)},1^{O(g)}) \to \mathbb{N}_0$$
be the number of closed saddle connections at a zero of order $m_i$ on $(X,\omega)$ with lengths in $[\frac{a}{\sqrt{g}},\frac{b}{\sqrt{g}}]$.
From Proposition \ref{order m} we have $L_{g,i,[a,b]}$ converges in distribution to a Poisson distribution of mean $\lambda_{i,[a,b]}$, where
$$\lambda_{i,[a,b]}=\frac{c_i(m_i^2-1)\pi}{2}(b^2-a^2).$$

Fix some $i,j\in \{1,2,\cdots,k-1\}$,
let
$$\hat{L}^{i,j}_{g,[a,b]}=(L_{g,i,[a,b]},L_{g,j,[a,b]}).$$
Similarly if we want to compute the distribution of $\hat{L}^{i,j}_{g,[a,b]}$, we need to consider the set 
$$(X,\omega,\Gamma_i,\Gamma_j),$$
where $(X,\omega)\in \mathcal{C}^{1,0}_{g}(m_1^{O_1(g)},\cdots,m_k^{O_k(g)},1^{O(g)})$, 
and $\Gamma_i$ is an ordered list of $n_i$ closed saddle connections at a zero of order $m_i$.
Again from \ref{exception}, these saddle connections can be collapsed simultaneously except for a measure zero set.
Then we have $L_{g,i,[a,b]}$ and $L_{g,j,[a,b]}$ converge to independent Poisson distributions of means $\lambda_{i,[a,b]}$ and $\lambda_{j,[a,b]}$.

Let $L_{g,[a,b]}$ be the number of closed saddle connections on $(X,\omega)$ with lengths in $[\frac{a}{\sqrt{g}},\frac{b}{\sqrt{g}}]$,
where $(X,\omega)\in \H(m_1^{O_1(g)},\cdots,m_k^{O_k(g)},1^{O(g)})$.
From the property of Poisson distribution,  we have $L_{g,[a,b]}=\sum\limits_{i=1}^{k}L_{g,i,[a,b]}$ converges in distribution to a Poisson distribution of mean $\lambda_{[a,b]}$, where
$$\lambda_{[a,b]}=\sum_{i=1}^{k} \lambda_{i,[a,b]}=\sum_{i=1}^{k-1}\frac{c_i}{2}(m_i^2-1)\pi(b^2-a^2).$$

Then we have proved Theorem \ref{general m}.

\bibliographystyle{plain}

\bibliography{bibliography}
\end{document}